\documentclass[12pt]{article}

\usepackage{fullpage}
\usepackage{graphicx}

\usepackage[affil-it]{authblk}

\usepackage{amsmath}
\usepackage{amstext, amssymb}
\usepackage{eufrak}
\usepackage[mathscr]{euscript}

\renewcommand{\d}[1]{\ensuremath{\operatorname{d}\!{#1}}}

\begin{document}

\title{Asymptotic Distribution of Centralized $r$ When Sampling from Cauchy}
\author[1]{Veson Lee}
\author[1]{Jan Vrbik}
\affil[1]{Department of Mathematics and Statistics\\ Brock University, Canada}
\date{September 15, 2018}
\maketitle

\begin{abstract}
Assume that $X$ and $Y$ are independent random variables, each having a Cauchy distribution with a known median. Taking a random independent sample of size $n$ of each $X$ and $Y$, one can then compute their centralized empirical correlation coefficient $r$. Analytically investigating the sampling distribution of this $r$ appears possible only in the large $n$ limit; this is what we have done in this article, deriving several new and interesting results.
\end{abstract}

\section{Introduction}
It can be easily shown, based on Central Limit Theorem, that the sampling distribution of the usual sample correlation coefficient further multiplied by $\sqrt{n}$ i.e.
\begin{equation*}
\frac{\sqrt{n}\cdot\sum_{i=1}^{n}\left(X_{i}Y_{i}-\frac{\sum_{i=1}^{n}X_{i}\cdot\sum_{i=1}^{n}Y_{i}}{n}\right)}{\sqrt{\sum_{i=1}^{n}X_{i}^{2}-\frac{\left(\sum_{i=1}^{n}X_{i}\right)^{2}}{n}}\sqrt{\sum_{i=1}^{n}Y_{i}^{2}-\frac{\left(\sum_{i=1}^{n}Y_{i}\right)^{2}}{n}}}
\end{equation*}
tends to, when $n\to\infty$ and $X_{i}$ and $Y_{i}$ are independent (both within and between), the standard Normal distribution whenever both $X$ and $Y$ have finite means and variances. The situation changes dramatically when sampling from a Cauchy distribution. Investigating what happens in that case is rather difficult; to simplify the task, we assume that the medians of $X$ and $Y$ are known, and can therefore be subtracted from the $X_{i}$ and $Y_{i}$ values. This amounts to assuming that the two Cauchy distributions have zero medians and it is then sufficient to define what we call the centralized sample correlation coefficient as 
\begin{equation}
r_{c}=\frac{\sum_{i=1}^{n}X_{i}Y_{i}}{\sqrt{\sum_{i=1}^{n}X_{i}^{2}}\sqrt{\sum_{i=1}^{n}Y_{i}^{2}}}
\label{eq:centered_r}
\end{equation}
It is obvious that the second parameter (the quartile deviation) of each of the two Cauchy distributions cancels out of the last expression; we can thus assume (without a loss of generality), that both the $X_{i}$ and $Y_{i}$ are drawn independently from a Cauchy distribution with median equal to 0 and the quartile deviation equal to 1. The objective of this article is to find the asymptotic (i.e. large-$n$) behaviour of $r_{c}$.  This centralized $r_{c}$ is also known to others in the computer science, big data and data science fields as the cosine similarity measure. It has applications in text mining, data mining and information retrieval $\text{\cite{key-1,key-2}}$; an investigation of the statistical distribution related to this centralized $r_{c}$ has been done by $\text{\cite{key-3}}$.

\section{Related Sampling Distributions}
To obtain its asymptotic distribution, we must first explore the distribution of the individual terms of $r_{c}$. We start by quoting the probability density function (PDF) and the characteristic function (CF) of each $X_{i}$ and $Y_{i}$: 
\begin{equation*}
f\left(x\right)=\frac{1}{\pi\left(x^{2}+1\right)}
\end{equation*}
and 
\begin{equation*}
\varphi\left(t\right)=\exp{\left(-\left|t\right|\right)}
\end{equation*}
respectively (these are well-known results $\text{\cite{key-4}}$).

This implies that the PDF each $X_{i}^{2}$ and $Y_{i}^{2}$ (denoted $Z$) is given by 
\begin{equation}
f\left(z\right)=\begin{cases}
\frac{\displaystyle 1}{\displaystyle\pi\sqrt{z}\left(z+1\right)} & z>0\\
0 & \textrm{elsewhere}
\end{cases}
\end{equation}
and the distribution function of $Z$ is given by 
\begin{equation*}
F_{X^{2}}\left(z\right)=\Pr\left(X^{2}\leq z\right)=\frac{1}{\pi}\int_{-\sqrt{z}}^{\sqrt{z}}\frac{\d x}{1+x^{2}}=\frac{2}{\pi}\arctan\sqrt{z}
\end{equation*}

The corresponding CF is then
\begin{equation}
\varphi_{X^{2}}\left(t\right)=\exp\left(-it\right)\left[1-\textrm{erf}{\left(\sqrt{-it}\right)}\right]
\label{eq:cf_cauchy_sq}
\end{equation}
and can be found by taking the appropriate Fourier transform of the PDF $\text{\cite{key-5}}$; here we rely on computer software (such as Maple or Mathematica) to provide these.

To obtain the asymptotic CF of each $\lim_{n\to\infty}\frac{\sum_{i=1}^{n}X_{i}^{2}}{n^{2}}$ and $\lim_{n\to\infty}\frac{\sum_{i=1}^{n}Y_{i}^{2}}{n^{2}}$ (denoted $W$), we have to raise $\eqref{eq:cf_cauchy_sq}$ to the power of $n$, replace $t$ by $\frac{t}{n^{2}}$ and then take the $n\to\infty$ limit of the resulting expression (note that only by dividing by $n^{2}$ can one reach a finite limit - that is how these `normalizing factors' are found in general).

This yields 
\begin{equation*}
\varphi_{W}\!\left(t\right)=\exp\left(-2\sqrt{\frac{-it}{\pi}}\right) 
\end{equation*}
since $\textrm{erf}{\left(x\right)}\simeq\frac{2x}{\sqrt{\pi}}$ for small $x$. The appropriate inverse Fourier transform converts this CF to the corresponding PDF, namely 
\begin{equation}
f\!\left(w\right)=\begin{cases}
\frac{\displaystyle\exp\left(-\frac{1}{\pi w}\right)}{\displaystyle\pi w^{\frac{3}{2}}} & w>0\\
0 & \textrm{elsewhere}
\end{cases}
\end{equation}
One can show via Monte Carlo simulation that this constitutes a fairly accurate approximation for the PDFs of $\frac{\sum_{i=1}^{n}X_{i}^{2}}{n^{2}}$ and $\frac{\sum_{i=1}^{n}Y_{i}^{2}}{n^{2}}$ even for relatively small $n$ values ($n\geq30$). This is due to the fact that the corresponding CF can be expanded in increasing powers of $\frac{1}{n}$. The error of the approximation is thus of the $O\left(\frac{1}{n}\right)$ type which is faster than $O\left(\frac{1}{\sqrt{n}}\right)$ of the Central Limit Theorem.

The PDF of W can be readily converted into the (asymptotic) PDF of each $\lim_{n\to\infty}\frac{n}{\sqrt{\sum_{i=1}^{n}X_{i}^{2}}}$ and $\lim_{n\to\infty}\frac{n}{\sqrt{\sum_{i=1}^{n}Y_{i}^{2}}}$ (denoted $U$ and equal to $\frac{1}{\sqrt{W}}$) using a simple univariate transformation, yielding 

\begin{equation}
f\!\left(u\right)=\begin{cases}
\frac{\displaystyle\exp{\left(-\frac{u^{2}}{\displaystyle\pi}\right)}}{\displaystyle\pi} & u>0\\
0 & \textrm{elsewhere}
\end{cases}
\label{eq:pdf_denom_asymp}
\end{equation}
resulting in a half-normal distribution.

Finally, each of the $X_{i}Y_{i}$ (denoted $S$) has a CF found by 
\begin{equation*}
\int_{-\infty}^{\infty}
\int_{-\infty}^{\infty}
\exp(ixyt)\cdot\frac{1}{\pi\left(x^{2}+1\right)}\cdot\frac{1}{\pi\left(y^{2}+1\right)}
\d x
\d y
\end{equation*}
resulting in 
\begin{equation*}
\varphi_{S}\left(t\right)=\begin{cases}
\frac{\displaystyle2\left(\sin{\left|t\right|}\textrm{Ci}{\left|t\right|}-\cos{\left|t\right|}\textrm{Si}{\left|t\right|}\right)}{\displaystyle\pi}+\cos{\left|t\right|} & t\neq0\\
1 & t=0
\end{cases}
\end{equation*}
which corresponds to the relatively simple PDF 
\begin{equation}
f\!\left(s\right)=
\begin{cases}
\frac{\displaystyle\ln\left(s^{2}\right)}{\displaystyle\pi^{2}\left(s^{2}-1\right)} & v=\mathbb{R}\setminus\left\{ -1,0,1\right\} \\
\frac{\displaystyle 1}{\displaystyle\pi^{2}} & v=\left\{ -1,1\right\} \\
0 & \textrm{elsewhere}\\
\end{cases}
\end{equation}

Assuming that each $X_{i}Y_{i}$ (and also their sum) are asymptotically independent of $\frac{n}{\sqrt{\sum_{i=1}^{n}X_{i}^{2}}}$ and $\frac{n}{\sqrt{\sum_{i=1}^{n}Y_{i}^{2}}}$ (something we have been able to verify only empirically), we now find the CF of 
\begin{equation}
X_{i}Y_{i}\cdot\frac{n}{\sqrt{\sum_{i=1}^{n}X_{i}^{2}}}\cdot\frac{n}{\sqrt{\sum_{i=1}^{n}Y_{i}^{2}}}\simeq S\cdot U_{1}\cdot U_{2}
\label{eq:rv_rc}
\end{equation}
where $U_{1}$ and $U_{2}$ are independent random variables, each having an asymptotic PDF of $\eqref{eq:pdf_denom_asymp}$. This resulting approximate CF is computed by 
\begin{equation*}
\frac{4}{\pi^{4}}
\int_{0}^{\infty}
\int_{0}^{\infty}
\int_{-\infty}^{\infty}
\cos\left(u_{1}u_{2}st\right)
\cdot\exp\!\left(-\frac{u_{1}^{2}+u_{2}^{2}}{\pi}\right)\cdot\frac{\ln{s^{2}}}{s^{2}-1} \d s \d u_{1} \d u_{2}
\end{equation*}
where the usual $\exp(iu_{1}u_{2}st)$ has beeen replaced by $\cos(u_{1}u_{2}st)$, since the resulting distribution is symmetric, implying that it is CF has no imaginary part.

To simplify the $\d u_{1} \d u_{2}$ integration, we perform it in the usual polar coordinates (denoted $R$ and $\Theta$) getting 
\begin{align}
& 
\frac{4}{\pi^{4}}
\int_{0}^{\infty}
int_{0}^{\frac{\pi}{2}}
\int_{-\infty}^{\infty}
R\cos\!\left(stR^{2}\frac{\sin{2\Theta}}{2}\right) \cdot
\exp\left(-\frac{R^{2}}{\pi}\right)\cdot\frac{\ln{s^{2}}}{s^{2}-1}
\d s \d \Theta \d R\nonumber \\
 & =
\frac{2}{\pi^{3}}
\int_{0}^{\frac{\pi}{2}}
\int_{-\infty}^{\infty}
\frac{1}{1+\left(\pi st\frac{\sin{2\Theta}}{2}\right)^{2}}
\cdot
\frac{\ln{s^{2}}}{s^{2}-1} \d s \d\Theta\nonumber \\
 & =
\frac{2}{\pi}
\int_{0}^{\frac{\pi}{2}}
\frac{1+\left|t\right|\sin{2\Theta}\ln{\left(\pi\left|t\right|\frac{\sin2\Theta}{2}\right)}}{1+\left(\pi t\frac{\sin{2\Theta}}{2}\right)^{2}}
\d\Theta
\label{eq:cf_rc_term}
\end{align}
\section{Asymptotic PDF of $r_{c}$}

To get an approximate CF of 
\begin{equation}
\frac{\sum_{i=1}^{n}X_{i}Y_{i}}{n\ln n}\cdot\frac{n}{\sqrt{\sum_{i=1}^{n}X_{i}^{2}}}\cdot\frac{n}{\sqrt{\sum_{i=1}^{n}Y_{i}^{2}}}. \label{eq:rv_rc_sum} 
\end{equation}
Note the unusual normalizing factor -- the only way to achieve a finite limit. We need to first raise $\eqref{eq:cf_rc_term}$ to the power of $n$ while simultaneously replacing $t$ by $\frac{t}{n\ln{n}}$, taking the $n\to\infty$ limit and finally evaluating the remaining integral. These operations can be carried out in any order in this particular case - something not true in general.

This yields before the $\d\Theta$ integration
\begin{equation}
\exp\left(-\left|t\right|\sin{2\Theta}\right) + \mathcal{O}\!\left(\frac{1}{\ln{n}}\right)
\label{eq:cf_rc_asymp}
\end{equation}
Since the error of this approximation is proportional to $\frac{1}{\ln{n}}$, the actual convergence of the subsequent result is expected to be rather slow (reaching a good accuracy only when $n$ is at least a few hundred). One can recover a part of this error (trying to recover the full $\frac{1}{\ln{n}}$ proportional term would make the resulting expression too cumbersome) by replacing $\eqref{eq:cf_rc_asymp}$ with the more accurate
\begin{equation*}
\exp\left(-a\left|t\right|\sin{2\Theta}\right)
\end{equation*}
where
\begin{equation*}
a=1-\frac{\ln{\pi}}{\ln{n}} .
\end{equation*}

To complete the computation, we evaluate 
\begin{equation*}
\varphi_{r_{c}}\left(t\right)=\frac{2}{\pi}\int_{0}^{\frac{\pi}{2}}\exp{\left(-a\left|t\right|\sin{2\Theta}\right)}\d\Theta=I_{0}{\left(a\left|t\right|\right)}-L_{0}{\left(a\left|t\right|\right)}
\end{equation*}
where $I_{0}$ and $L_{0}$ are the Bessel and modified Struve functions respectively. Converting to PDF results in: 
\begin{equation}
f\!\left(r_{c}\right)=\frac{2\ln\left(a+\sqrt{a^{2}+r_{c}^{2}}\right)-\ln r_{c}^{2}}{\pi^{2}\sqrt{a^{2}+r_{c}^{2}}}
\label{eq:pdf_rc_aymp}
\end{equation}
which is our final answer for the approximate distribution of the sum of $\eqref{eq:rv_rc_sum}$ terms, which is the same as
\begin{equation}
\frac{n\cdot r_{c}}{\ln{n}}
\label{eq:rv_rc_normalized} .
\end{equation}
Note that similarly to the Cauchy distribution itself $\eqref{eq:pdf_rc_aymp}$ has an indefinite mean and an infinite variance.

\section{Monte Carlo Verification}
We will now verify the accuracy of our approximation by randomly generating $100,000$ values of $\eqref{eq:rv_rc_normalized}$, using $n=400$. The following Mathematica program does this and plots the corresponding histogram together with our approximate PDF $\eqref{eq:pdf_rc_aymp}$.
\begin{verbatim}
G[x_, y_] := x.y/Sqrt[x.y y.y]
n = 400;
superN = 100 000;
data = {};
Do[SeedRandom[];
  dataX = RandomReal[CauchyDistribution[], n];
  dataY = RandomReal[CauchyDistribution[], n];
  AppendTo[data, (n/Log[n]) * Apply[G, {dataX, dataY}]], {superN} ]
  
a = 1 - Log[Pi]/Log[n];
pdf = (2*Log[a + Sqrt[a^2 + x^2]] - Log[x^2]) / (Pi^2*Sqrt[a^2 + x^2]);
Show[Histogram[data, {-4, 4, 0.25}],
  Plot[0.25*superN*pdf, {x,-4,-4}, PlotRange -> {0,0.2*superN}]]
\end{verbatim}

\begin{figure}[h!]\centering
\includegraphics[width=.7\textwidth]{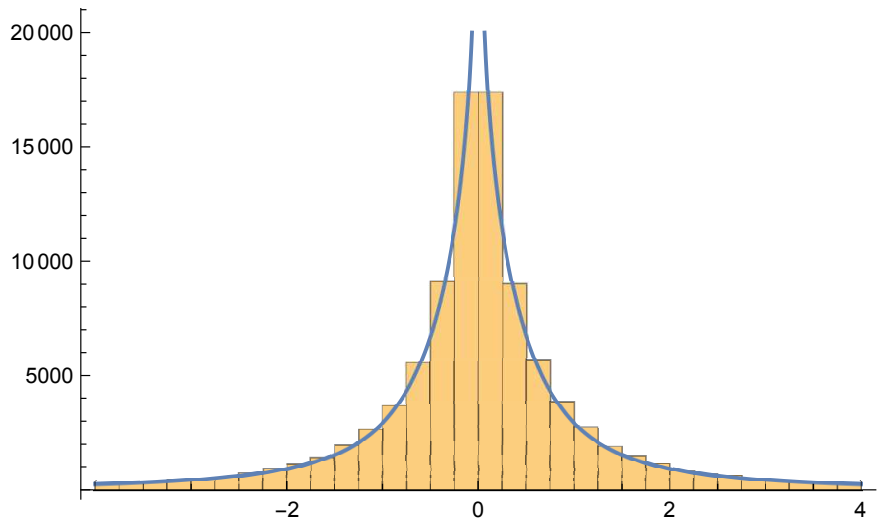}
\caption{Histogram vs.\ Approximate PDF.}
\label{thisone}
\end{figure}

The two results are in good agreement, as can be seen in the Figure \ref{thisone}.

\section{Conclusion}
We have derived an asymptotic distribution of the centralized sample correlation coefficient when sampling from Cauchy distribution, discovering some of its unusual properties. As a byproduct, several other interesting distributions have been introduced in the process.

\bibliographystyle{plain}

\end{document}